\documentclass{amsart}
\usepackage{amssymb,latexsym}

\numberwithin{equation}{section}

\newtheorem{theorem}{Theorem}[section]
\newtheorem{proposition}[theorem]{Proposition}
\newtheorem{corollary}[theorem]{Corollary}

\newtheorem{lemma}[theorem]{Lemma}

\def\mathbb\bf
\def\mathbf\bf
\def\mathcal\cal
\def\CC{{\bf C}}

\def\ZZ{{\bf Z}}
\def\gl{{\rm GL}}
\def\sp{{\rm Sp}}
\def\FF{{\mathcal F}}
\def\CFF{\sp\FF}

\def\fl{{\rm Fl}}
\def\cfl{\sp\fl}

\def\<{\langle}
\def\>{\rangle}
\def\opp{{\rm op}}
\def\sym{{\rm sym}}
\def\id{{\rm id}}
\def\cpr{{\rm cpr}}

\newcommand{\choos}[2]{\!
\mbox{\footnotesize $
{\left(\!\!\! \begin{array}{c} #1 \\ #2 \end{array}\!\!\!\right)}
$ }\!\!}

\def\aa{{\mathbf{a}}}
\def\bb{{\mathbf{b}}}
\def\cc{{\mathbf{c}}}
\def\dd{{\mathbf{d}}}
\def\ee{{\mathbf{e}}}
\def\ff{{\mathbf{f}}}
\def\kk{{\mathbf{k}}}

\def\Lam{\Lambda}

\begin{document}

\title{Symplectic multiple flag varieties of finite type}

\author{Peter Magyar}

\address{Dept. of Mathematics,
Northeastern University\\
Boston, Massachusetts 02115, USA}
\email{pmagyar@lynx.neu.edu}

\author{Jerzy Weyman}

\email{weyman@neu.edu}

\author{Andrei Zelevinsky}
\email{andrei@neu.edu}

\thanks{The research of Peter Magyar,
Jerzy Weyman and Andrei Zelevinsky
was supported in part by an NSF Postdoctoral Fellowship and
NSF grants \#DMS-9700884 and
\#DMS-9625511, respectively.}

\subjclass{14M15, 16G20, 14L30.}

\date{\today }

\maketitle

\vspace{-1em}
\noindent
\centerline{\it Dedicated to Professor David Buchsbaum on his retirement.}
\\[1em]
In our paper \cite{mwz}, we examined
the following question in the case of the group $G = \gl_m$:
\\[.5em]
{\bf Problem:} Given a reductive algebraic group $G$,
find all $k$-tuples of parabolic subgroups
$(P_1,\ldots,P_k)$ such that the product of flag varieties
$G/P_1 \times\cdots \times G/P_k$
has finitely many orbits under the diagonal action of $G$.
In this case we call $G/P_1 \times\cdots \times G/P_k$ a
{\it multiple flag variety of finite type}.
\\[.5em]
In this paper, we solve this problem for
the symplectic group $G = \sp_{2n}$.
We also give a complete enumeration of the orbits, and
explicit representatives for them.
The cases in our classification where
one of the parabolics is a Borel subgroup, $P_1=B$,
are exactly those for which $G/P_2\times \cdots \times G/P_k$ is a
spherical variety under the diagonal action of $G=\sp_{2n}$,
and our results specialize to classify the $B$-orbits on these
spherical varieties.

Our main tool is, as in \cite{mwz},
the algebraic theory of quiver representations.
Rather unexpectedly, it turns out that two multiple 
symplectic flags lie in the same $\sp_{2n}$-orbit if and only if
they lie in the same $\gl_{2n}$-orbit
(a consequence of a general result in
Proposition~\ref{pr:abstract orbit correspondence} below).
This allows us to reduce our problem for $\sp_{2n}$
to one about $\gl_{2n}$, which we solve using the results and
the quiver techniques of \cite{mwz}.

All of our methods extend in an obvious way to
the orthogonal groups $G = {\rm SO}_m$,
but the combinatorics of orbits becomes much more complicated and
we will not present it here.
Our work intersects with that of Littelmann \cite{littelmann},
who solved our problem for an arbitrary $G$, but with
the restrictions that $k = 3$, $P_1 = B$, and $P_2$, $P_3$
are maximal parabolic subgroups.

\textsc{Acknowledgments.}
Some crucial technical results in our paper are influenced
by a work in progress \cite{dw}.
We thank Harm Derksen for generously sharing his insights with us,
and Michel Brion, Dominique Luna, and Dmitry Panyushev for
their helpful interest in our work.
\\[1em]
{\small
\textsc{contents}
{\bf 1.} Main results \quad
{\bf 1.1.} Classification of finite types \quad
{\bf 1.2}  Flag categories \quad
{\bf 1.3}  Classification of orbits \quad
{\bf 1.4}  Enumeration of orbits \qquad
{\bf 2.} Proofs \quad
{\bf 2.1.} Proof of Theorem 1.3(i) \quad
{\bf 2.2.} Proof of Theorem 1.3(ii) \quad
{\bf 2.3.} Proof of Theorems 1.1 and 1.2 \quad
{\bf 2.4.} Proof of Theorem 1.4
}

\section{Main results}

\subsection{Classification of finite types}
\label{sec:classification}

We recall the notation of \cite{mwz}.
A \emph{composition}
of a positive integer $m$ is a
sequence of non-negative integers
$\aa = (a_1, \ldots, a_p)$ whose sum is equal to $m$.
The
components $a_i$ are called \emph{parts} of $\aa$.
Define the
\emph{opposite} of $\aa = (a_1, \ldots, a_p)$,
by $\aa^{\opp}=(a_p, \ldots, a_1)$, and
say that $\aa$ is {\it symmetric} if $\aa = \aa^{\opp}$.
Let $(a^k)$ denote the composition with $k$ parts all equal to $a$.

Thoughout this paper, all vector spaces are
over a fixed algebraically closed field.
For a vector space $V$ of dimension $m$ and a composition $\aa$ of $m$,
we denote by $\fl_{\aa} (V)$ the variety of flags
$A = (0 = A_0 \subset A_1 \subset \cdots \subset A_p = V)$
of vector subspaces in $V$ such that
$$
\dim (A_{i}/A_{i-1}) = a_i \quad (i = 1, \ldots, p) \ .
$$
Now let $V$ be a $2n$-dimensional symplectic vector space
possessing a non-degenerate alternating bilinear form $\<\,,\,\>$.
The group of automorphisms of $V$ preserving the form is $\sp(V) = \sp_{2n}$.
A subspace $U\subset V$ is \emph{isotropic} if $\<U,U\>=0$.
For a symmetric composition $\aa$ of $2n$, we denote by
$\cfl_{\aa} (V)$ the variety of flags of dimension vector $\aa$
in $V$ that are formed by isotropic subspaces and their orthogonals:
$$
\cfl_{\aa}(V) =
\{ A \in \fl_{\aa}(V) \mid \<A_i, A_{p-i}\> = 0\ \mbox{for all}\ i\}.
$$
This is a standard realization of a partial flag variety $\sp(V)/P$.
(See \cite[\S23.3]{fulton-harris}.)
The complete flag variety
$\sp(V)/B$ corresponds to the compositon $\aa=(1^{2n})$.

A tuple of symmetric compositions $\dd=(\aa_1, \ldots, \aa_k)$ of
the same number $2n$ is said to be \emph{of symplectic finite type}
($\sp$-finite)
if the group $\sp (V)$, acting diagonally, has finitely many orbits in
the \emph{multiple flag variety}
$\cfl_{\dd}(V)=\cfl_{\aa_1} (V) \times \cdots \times \cfl_{\aa_k} (V)$.
We will classify all such tuples.

We say that a composition is \emph{trivial}
if it has only one part, $\aa=(2n)$.
Then the corresponding flag variety $\cfl_{\aa}$
consists of a single point,
so adjoining any number of trivial compositions to a tuple
gives essentially the same multiple flag variety
and does not affect the finite-type property.

\begin{theorem}
\label{th:k=3}
If a tuple of non-trivial symmetric compositions $(\aa_1, \ldots, \aa_k)$ is
of symplectic finite type then $k \leq 3$.
\end{theorem}

Thus we only need to classify \emph{triples} of symplectic finite type.
We will write $\dd=(\aa, \bb, \cc)$ instead of
$(\aa_1, \aa_2, \aa_3)$.

The vanishing of some part $a_i=0$ of $\aa$
means that in any flag $A \in {\fl}_{\aa}$, the subspace
$A_i$ coincides with $A_{i-1}$.  Thus, removing zero parts
$a_i$ and $a_{p+1-i}$ from a symmetric composition $\aa$
does not change $\fl_{\aa}$ or $\cfl_{\aa}$ up to isomorphism.
Given $(\aa,\bb,\cc)$, let $p, q$ and $r$ denote their respective
numbers of non-zero parts.
Assume without loss of generality that $p \leq q \leq r$.

\begin{theorem}
\label{th:finite-type-list}
A triple $(\aa, \bb, \cc)$ of symmetric compositions of $2n$
is of symplectic finite type if and only if it belongs to
one of the following classes:
\\
$(\sp A_{q,r})$\ : \,\, $(p,q,r) = (1,q,r)$, $1 \leq q \leq r$.
\\[.5em]
$(\sp D_{r+2})$\ : \,\, $(p,q,r) = (2,2,r)$, $2 \leq r$.
\\[.5em]
$(\sp E_6)$\  : \,\, $(p,q,r) = (2,3,3)$.
\\[.5em]
$(\sp E_7)$\ : \,\, $(p,q,r) = (2,3,4)$.
\\[.5em]
$(\sp E_8)$\ : \,\, $(p,q,r) = (2,3,5)$.
\\[.5em]
$(\sp E^{(b)}_{r+3})$\ : \,\, $(p,q,r) = (2,3,r)$, $3 \leq r$,
$\bb$  {\rm has non-zero  parts} $(1, 2n-2, 1)$.
\\[.5em]
$(\sp Y_{r+4})$\ : \,\, $(p,q,r) = (3,3,r)$, $3 \leq r$,
$\aa, \bb$, {\rm or} $\cc$
{\rm has non-zero  parts} $(1, 2n-2, 1)$.
\end{theorem}

The labels are taken from \cite{mwz},
and refer to Dynkin graphs associated to the first five cases.
Note that except for type ($\sp Y$), the dimension vectors $\dd$
of all the above types also appear on our list of $\gl$-finite dimensions
in \cite[Theorem 2.2]{mwz}.  That is, not only does $\sp_{2n}$ have
finitely many
orbits on $\cfl_{\dd}$, but $\gl_{2n}$ has finitely many orbits on $\fl_{\dd}$.

As in \cite{mwz}, the type ($\sp A$) covers all symplectic
multiple flag varieties with only one or two non-trivial factors.
Note that if $p = 2$ then $\aa = (n,n)$, and the corresponding variety
$\cfl_{(n,n)} (V)$ is the variety of all Lagrangian subspaces in $V$.
Thus the case ($\sp D$) covers triple symplectic flag varieties
in which two of the factors are Lagrangian Grassmannians, and the third
factor is arbitrary.

In general, each flag in a symplectic flag variety is completely
determined by its ``lower half" consisting of isotropic subspaces.
Thus, the cases ($\sp E_6$), ($\sp E_7$), and ($\sp E_8$) correspond to
triple flag varieties in which:
the first flag contains a single Lagrangian subspace;
the second flag contains a single isotropic subspace;
and the third flag contains at most two isotropic subspaces.
The cases ($\sp E^{(b)}_{r+3}$) and ($\sp Y_{r+4}$) correspond to
triple flag varieties in which: the first flag contains
a single isotropic subspace; the second flag
contains a single line (automatically isotropic);
and the third flag contains arbitrary isotropic subspaces.

Recall that the variety $\cfl_{(\aa,\bb)}$ is {\it spherical}
whenever $(\aa, \bb, (1^{2n}))$ is of $\sp$-finite type.
For most of the above $\sp$-finite triples $(\aa,\bb,\cc)$,
the triple $(\aa,\bb,(1^{2n}))$ is also $\sp$-finite,
and the variety $\cfl_{(\aa,\bb)}$ is spherical.
The exceptions are the cases ($\sp E_6$), ($\sp E_7$), and
($\sp E_8$), provided $2n \geq 6$ and $(\aa,\bb,(1^{2n}))$ is not of type 
($\sp E^{(b)}_{r+3}$).  
These latter cases go beyond the
scope of Littelmann's classification \cite{littelmann}.

\subsection{Flag categories}

For each $\sp$-finite triple $\dd=(\aa,\bb,\cc)$,
we describe the $\sp(V)$-orbits on the triple flag variety
$\cfl_{\dd}(V)$.
Remarkably, we can do so in the same categorical framework
as in the case of $G=\gl_{m}$.
Later, in Corollary \ref{co:combinatorial-orbit-param},
we give the parametrization of the orbits
in purely combinatorial terms.

For a composition $\aa = (a_1, \ldots, a_p)$, we write
$
|\aa| = a_1 + \cdots + a_p.
$
The number $p$ of parts of $\aa$ will be denoted $\ell (\aa)$,
called the \emph{length} of $\aa$.
For any positive integers $p, q$, and $r$, we consider an
additive semigroup of triples of compositions:
$$
\Lam_{pqr} = \{(\aa, \bb, \cc) \mid
(\ell (\aa), \ell (\bb), \ell (\cc)) = (p,q,r) \mbox{ and }
|\aa| = |\bb| = |\cc|\}.
$$
(Here, in contrast to the notation of Theorem~\ref{th:finite-type-list},
the numbers $p,q,r$ include the zero parts of $\aa$, $\bb$, $\cc$.)

Introduce the additive category $\FF=\FF_{pqr}$ whose
objects are families $(V; A, B, C)$,
where $V$ is any finite-dimensional vector space, and
$(A,B,C)$ is a triple of flags in $V$ belonging to any
$\fl_{\aa} (V) \times \fl_{\bb} (V)  \times \fl_{\cc} (V)$
with $(\aa, \bb, \cc) \in \Lam_{pqr}$.
The triple $\dd = (\aa, \bb, \cc)$ is called
the \emph{dimension vector} of $(V; A, B, C)$.
A \emph{morphism} in $\FF$ from $(V; A, B, C)$ to $(V'; A', B', C')$
is a linear map $f: V \to V'$ such that
$f(A_i) \subset A'_i, \, f(B_i) \subset B'_i$, and $f(C_i) \subset C'_i$
for all $i$.
Direct sum of objects is taken componentwise
on each member of each flag.

Each triple flag in
$\fl_{\aa} (V) \times \fl_{\bb} (V)  \times \fl_{\cc} (V)$
corresponds naturally to an object of $\FF$, and
$\gl(V)$-orbits in the triple flag variety
are naturally identified with isomorphism classes of objects
in $\FF$ with dimension vector $(\aa, \bb, \cc)$.
The advantage of translating the $\gl$-orbit problem
into the additive category is that each object
splits \emph{uniquely} into indecomposable objects,
so that an isomorphism class is uniquely specified by the
multiplicities of its indecomposable summands (cf.~\cite{mwz}).
Thus to classify all isomorphism classes in a given dimension it
is enough to find the \emph{indecomposable} isomorphism classes in all
smaller dimensions.  But indecomposables are rather rare, and
this becomes a tractable and familiar problem in the theory of quivers.

Next we translate the $\sp$-orbit problem into categorical terms.
Let $\CFF=\CFF_{pqr}$ be the full subcategory of $\FF_{pqr}$ consisting
of the objects $(V;A,B,C)$ which have symmetric dimension vector and
such that $V$ admits a non-degenerate symplectic form $\<\,,\,\>$ with
$\<A_{p-i},A_i\>=\<B_{q-i},B_i\>=\<C_{r-i},C_i\>=0$ for all $i$.
Now fix a symplectic form on $V$.
Clearly a triple flag in $\cfl_{\aa}(V) \times \cfl_{\bb}(V)
\times \cfl_{\cc}(V)$ may be thought of as an object of $\CFF_{pqr}$,
and if two such triple flags are in the same $\sp(V)$-orbit, then
they are equivalent as objects of $\CFF_{pqr}$.
Also, since all symplectic forms on $V$ are conjugate, every
isomorphism class of $\CFF_{pqr}$ contains at least one $\sp(V)$-orbit
of our triple symplectic flag variety.
Our first key technical result is that each $\CFF$-class contains
\emph{exactly} one $\sp(V)$-orbit.

Our second key fact is a nice description of the indecomposable objects of
$\CFF$ in terms of indecomposables in $\FF$.
To give this description, we define a contravariant duality functor
$\ast$ on $\FF_{pqr}$.
For a single flag
$A = (A_1\subset \cdots \subset A_{p-1} \subset V)$, define
$A^* = ((V/A_{p-1})^* \subset \cdots \subset (V/A_{1})^*\subset V^*)$,
where $V^*$ denotes the dual vector space and
$(V/A_i)^*$ the subspace of linear forms vanishing on $A_i$.
Let $(V;A,B,C)^*=(V^*;A^*, B^*, C^*)$.
If an object has dimension
vector $\dd =(\aa,\bb,\cc)$, its dual has dimension vector
$\dd^{\opp} = (\aa^{\opp}, \bb^{\opp}, \cc^{\opp})$.

Notice that an object in $\CFF$ must be isomorphic to its dual
in $\FF$, since the symplectic form identifies
$V\stackrel{\sim}{\to}V^*$.
Clearly, not all self-dual objects of $\FF$ lie in $\CFF$ (for example,
$V$ might have odd dimension).
But it is easy to see that $I \oplus I^*$
is in $\CFF$ for any object $I$ of $\FF$;
indeed, if $F$ is any flag in $V$ then $F \oplus F^*$ is a symplectic
flag in $V \oplus V^*$ with respect to the symplectic form
having both $V$ and $V^*$ isotropic and inducing the natural pairing
between $V$ and $V^*$.

\begin{theorem}
\label{th:GL-Sp relationship}
(i) Two triple symplectic flags in $\cfl_{\dd}(V)$
lie in the same $\sp(V)$-orbit if and only if they
are isomorphic as objects in $\CFF_{pqr}$.\\
(ii) Let $J$ be an object of $\FF_{pqr}$.
Then $J$ is an indecomposable object of $\CFF_{pqr}$
if and only if \\
\mbox{} \quad (a) $J$ lies in $\CFF_{pqr}$ and $J$ is indecomposable in
$\FF_{pqr}$, or \\
\mbox{} \quad (b) $J \cong I\oplus I^*$ for some $\FF_{pqr}$-indecomposable $I$
not belonging to $\CFF_{pqr}$.
\end{theorem}

As a consequence of part (ii),  any object of $\CFF$ can be
\emph{uniquely} written as a sum of indecomposable objects of
$\CFF$.
Thus, as in the $\gl$ case, to classify orbits on a
multiple $\sp$-flag variety it is enough 
to find the $\CFF$-indecomposables which can appear as summands
of an object with the given dimension vector.

\subsection{Classification of orbits}

For a dimension vector $\dd = (\aa,\bb,\cc)$
with lengths $(p,q,r)$, we have seen
the $\sp$-orbits in $\cfl_{\dd}$
are naturally identified with direct sums
$\bigoplus_J m_J J$, where $J$ runs over all indecomposable
isomorphism classes in $\CFF_{pqr}$,
and the $m_J$ are non-negative integers with $\sum_J m_J \dim(J) = \dd$.
We proceed to exhibit representatives for
all $J$ which can appear in such a decomposition
for the dimensions $\dd$ of Theorem \ref{th:finite-type-list}.
This will give a classification of orbits, as well as explicit
triples of flags lying in each orbit.

We say a composition is \emph{compressed} if it has no zero parts,
and a flag with no repeated subspaces is also called compressed.
Let $\aa_{\rm cpr}$ denote the composition obtained from $\aa$
by removing all zero parts; let $A_{\rm cpr}$ denote the flag
obtained from $A$ by removing all repetitions of subspaces;
and let similar notation hold 
for tuples of compositions and multiple flags.
Even if an object in $\CFF_{pqr}$ is compressed,
its direct summand $J$ might not be compressed.
However, if we have a representative for $J_{\rm cpr}$
and we know $\dim(J)$, then we can immediately construct
a representative for $J$.  Thus in our list we need only
produce representatives for compressed $\CFF$-indecomposables $J=J_{\cpr}$.

As a further normalization, if an object $J = (V;A,B,C)$
has dimension $\dd=(\aa,\bb,\cc)$ where
$p\leq q\leq r$ does \emph{not} hold,
we switch $A$, $B$, and $C$ to get an object $J^{<}$
of dimension $\dd^<$ for which
the lengths of the flags \emph{are} non-decreasing.
(If there is more than one way to do such switching, choose an
arbitrary one.)

We adopt the following notation for $\CFF$-indecomposables $J$:\\
(1) Suppose that in dimension $\dd$ there are exactly $k$
isomorphism classes of $\CFF$-indecomposables which are
also indecomposable in $\FF$.  Then we denote these indecomposables
by $J=I_{\dd}^1, \ldots, I_{\dd}^k$ (or simply $J=I_{\dd}$ if $k=1$).\\
(2)  Suppose $\dd = \ee+\ee^{\opp}$,
and in dimension $\ee$ there is a unique isomorphism class of
$\FF$-indecomposables $I_{\ee}$ which is not in $\CFF$.  Then
$J=I_{\ee}\oplus I_{\ee}^*$ is a $\CFF$-indecomposable of dimension $\dd$,
which we denote $J = I_{\ee}^{\sym}$.

\begin{theorem}
\label{th:standard form}
An object $J$ is a $\CFF$-indecomposable summand of a $\CFF$ object
with symplectic finite dimension if and only if $J_{\rm cpr}^{<}$ is
isomorphic to
one of the objects in the table below.
(In the right-hand column of the
table, we abbreviate $I_{\dd}^j$ to $\dd^j$ and
$I_{\ee}^{\sym}$ to $\ee^{\sym}$, and we omit commas.)
$$
\begin{array}{cl}
\mbox{\rm dimensions}\ \dd & \mbox{\rm $\sp$-indecomposables}\ J
\\[.5em]
((2)(2)(2)) & ((1)(1)(1))^{\sym}
\\[.5em]
((2)(2)(11)) & ((1)(1)(10))^{\sym}
\\[.5em]
((2)(11)(11)) & ((1)(10)(10))^{\sym},\ ((1)(10)(01))^{\sym}
\\[.5em]
((11)(11)(11)) &
((11)(11)(11)),\\&
((10)(10)(10))^{\sym},\ ((10)(10)(01))^{\sym},\\&
((10)(01)(10))^{\sym},\ ((10)(01)(01))^{\sym}
\\[1em]
((22)(22)(121)) &
((11)(11)(110))^{\sym}
\\[.5em]
((22)(22)(1^4)) &
((11)(11)(1100))^{\sym},\
((11)(11)(1010))^{\sym}
\\[.5em]
((22)(121)(121)) &
((11)(110)(110))^{\sym},\
((11)(110)(011))^{\sym}
\\[.5em]
((22)(121)(1^4)) &
((22)(121)(1^4)),\\&
((11)(110)(1100))^{\sym},\
((11)(110)(1010))^{\sym},\\&
((11)(110)(0101))^{\sym},\
((11)(110)(0011))^{\sym}
\\[.5em]
((121)(121)(121)) &
((121)(121)(121)),\\&
((110)(110)(110))^{\sym},\
((110)(110)(011))^{\sym},\\&
((110)(011)(110))^{\sym},\
((110)(011)(011))^{\sym}
\\[.5em]
((121)(121)(1^4))&
((121)(121)(1^4))^1,\
((121)(121)(1^4))^2,\\&
((110)(110)(1100))^{\sym},\
((110)(110)(1010))^{\sym},\\&
((110)(110)(0101))^{\sym},\
((110)(110)(0011))^{\sym},\\&
((110)(011)(1100))^{\sym},\
((110)(011)(1010))^{\sym},\\&
((110)(011)(0101))^{\sym},\
((110)(011)(0011))^{\sym}
\end{array}  $$ $$ \begin{array}{cl}
((33)(2^3)(2^3))&
((21)(1^3)(1^3))^{\sym}
\\[.5em]
((33)(2^3)(2112))&
((21)(1^3)(1101))^{\sym},\
((21)(1^3)(1011))^{\sym}
\\[.5em]
((33)(2^3)(1221))&
((21)(1^3)(1110))^{\sym},\
((21)(1^3)(0111))^{\sym}
\\[.5em]
((33)(2^3)(11211))&
((33)(2^3)(11211)),\\&
((21)(1^3)(11100))^{\sym},\
((21)(1^3)(10110))^{\sym},\\&
((21)(1^3)(01101))^{\sym},\
((21)(1^3)(00111))^{\sym}
\\[.5em]
((141)(2^3)(2^3))&
((120)(1^3)(1^3))^{\sym}
\\[.5em]
((141)(2^3)(2112))&
((120)(1^3)(1101))^{\sym}\,
((120)(1^3)(1011))^{\sym}\,
\\[.5em]
((141)(2^3)(1221))&
((120)(1^3)(1110))^{\sym}\,
((120)(1^3)(0111))^{\sym}\,
\\[.5em]
((141)(2^3)(11211))&
((141)(2^3)(11211)),\\&
((120)(1^3)(11100))^{\sym},\
((120)(1^3)(10110))^{\sym},\\&
((120)(1^3)(01101))^{\sym},\
((120)(1^3)(00111))^{\sym}
\\[.5em]
((141)(2^3)(1^6))&
((141)(2^3)(1^6))^1,\
((141)(2^3)(1^6))^2,\\&
((120)(1^3)(111000))^{\sym},\
((120)(1^3)(110100))^{\sym},\\&
((120)(1^3)(101010))^{\sym},\
((120)(1^3)(100110))^{\sym},\\&
((120)(1^3)(011001))^{\sym},\
((120)(1^3)(010101))^{\sym},\\&
((120)(1^3)(001011))^{\sym},\
((120)(1^3)(000111))^{\sym}
\\[1em]
((44)(323)(2^4))&
((22)(211)(1^4))^{\sym}
\\[.5em]
((44)(323)(21212))&
((22)(211)(11101))^{\sym},\
((22)(211)(10111))^{\sym}
\\[.5em]
((44)(323)(12221))&
((22)(211)(11110))^{\sym},\
((22)(211)(01111))^{\sym}
\\[.5em]
((44)(242)(21212))&
((22)(121)(11101))^{\sym}
\\[.5em]
((44)(242)(12221))&
((22)(121)(11110))^{\sym}
\\[1em]
((55)(424)(2^5))&
((32)(212)(1^5))^{\sym}
\\[.5em]
((55)(343)(2^5))&
((32)(221)(1^5))^{\sym},\
((32)(122)(1^5))^{\sym}
\\[1em]
((66)(4^3)(32223))&
((3^2)(2^3)(21111))^{\sym}
\\[.5em]
((66)(4^3)(23232)&
((3^2)(2^3)(12111))^{\sym}
\end{array}
$$
\end{theorem}

Note that for all the above $\CFF$-indecomposables $J=(V;A,B,C)$,
the dimension of the ambient space $V$ is at most 12.
This is in contrast to the $\FF$-indecomposables of $\gl$-finite
type \cite[Theorem 2.9]{mwz}, which occur in all dimensions.

In all the above cases with $J=I_{\ee}^{\sym}$,
the dimension vector $\ee$ is $\gl$-finite
(i.e.~there are only finitely many $\gl$-orbits
in $\fl_{\ee}$).  Thus
we may read off an explicit representative for
$I_{\ee}$ from the list of $\FF$-indecomposables
of $\gl$-finite type in \cite[Theorem 2.9]{mwz}.
For the remaining indecomposable classes $J = I_{\dd}^i$
we give symplectic representatives below.

We present each indecomposable
$J=(V;A,B,C)$ with dimension vector $\dd=(\aa,\bb,\cc)$ as
follows.  The space $V=V_{2n}$ has basis $e_1,\ldots,e_{2n}$
(where $2n = |\aa|=|\bb|=|\cc|$),
with the standard symplectic form $\<e_i,e_{2n+1-i}\>
=-\<e_{2n+1-i},e_i\>=1$ for $i=1,\ldots,n$ and
$\<e_i,e_j\>=0$ otherwise.
We list explicit bases for the isotropic spaces in
$A$, $B$, and $C$ (the rest of the spaces being the orthogonals of
those given).
$$
\begin{array}{l@{\!}l}
I_{((11)(11)(11))}&\ =\ (V_2;\, (e_1+e_2),\ (e_1), \ (e_2)) \\[.4em]
I_{((22)(121)(1^4))}&\ =\ (V_4;\, (e_1+e_3+e_4, e_2 + e_4),\
(e_1), \ (e_4) \subset (e_3,e_4)) \\[.2em]
I_{((121)(121)(121))}&\ =\ (V_4;\, (e_1+e_2+e_4),\ (e_1), \ (e_4)) \\[.4em]
I_{((121)(121)(1^4))}^1&\ =\ (V_4;\, (e_1+e_2+e_4),\ (e_1), \ (e_4) \subset
(e_2, e_4)) \\[.4em]
I_{((121)(121)(1^4))}^2&\ =\ (V_4;\, (e_1+e_2+e_4),\, (e_1), \ (e_4) \subset
(e_3, e_4)) \\[.4em]
I_{((33)(2^3)(11211))}&\ =\ (V_6;\, (e_1+e_5+e_6, e_2+e_4+e_6, e_3+e_5),\
(e_1, e_2),	\ (e_6) \subset (e_5,e_6)) \\[.4em]
I_{((141)(2^3)(11211))}&\ =\ (V_6;\, (e_1+e_3+e_5+e_6),\ (e_1, e_2), \
(e_6) \subset (e_5,e_6)) \\[.4em]
I_{((141)(2^3)(1^6))}^1&\ =\ (V_6;\, (e_1+e_3+e_5+e_6), \ (e_1, e_2), \
(e_6) \subset (e_5,e_6) \subset (e_3, e_5,e_6))  \\[.4em]
I_{((141)(2^3)(1^6))}^2&\ =\ (V_6;\, (e_1+e_3+e_5+e_6), \ (e_1, e_2), \
(e_6) \subset (e_5,e_6) \subset (e_4, e_5,e_6)) \ .
\end{array}
$$

\subsection{Enumeration of orbits}

 From Theorem \ref{th:standard form} we may deduce the following
enumeration of orbits on multiple symplectic
flag varieties, similar to the Kostant partition function.

Let $\sp\Pi_{pqr}$ be the set of the symmetric triples
$\dd\in \Lam_{pqr}$ such that $\dd_{\cpr}^<$ is one of the
dimension vectors in the left-hand column of the table
in Theorem \ref{th:standard form}.
For each $\dd \in \sp\Pi_{pqr}$, let
$\mu_{\dd}$ be the number of $\sp$-indecomposables
of dimension $\dd_{\cpr}^<$ listed in the right-hand column next to
$\dd_{\cpr}^<$.
Define $\widetilde{\sp\Pi}_{pqr} \subset \sp\Pi_{pqr} \times \ZZ_{+}$ by
$$
\widetilde{\sp\Pi}_{pqr}
= \bigcup_{\dd \in \sp\Pi_{pqr}} \{\dd\} \times \{1,2,\ldots,\mu_{\dd}\}\ .
$$

\begin{corollary}
\label{co:combinatorial-orbit-param}
Let $\dd = (\aa,\bb,\cc) \in \Lam_{pqr}$ be a triple of
symmetric compositions.\\
(i) The $\sp_{2n}$-orbits on $\cfl_{\dd}$ are in natural bijection with
families of non-negative integers $(m_{\ee,i})$
indexed by $(\ee,i)\in \widetilde{\sp\Pi}_{pqr}$
such that
$$
\sum_{(\ee,i) \in \widetilde{\sp\Pi}} m_{\ee,i}\  \ee = \dd\ .
$$
(ii) The number of $\sp_{2n}$-orbits on $\cfl_{\dd}$ is:
$$
\sum_{(m_{\ee})}\, \prod_{\ee\in \sp\Pi} \choos{\mu_{\ee}+m_{\ee}-1}{
\mu_{\ee}-1},
$$
where the sum runs over all families of non-negative
integers $(m_{\ee})$ indexed by
$\ee\in \sp\Pi_{pqr}$ such that
$$
\sum_{\ee\in \sp\Pi} m_{\ee}\,\ee = \dd\ .
$$
\end{corollary}

Part (i) is a consequence of Theorem \ref{th:standard form}.
Part (ii) follows from (i), since the binomial coefficient
$\choos{\mu + m -1}{ \mu -1}$
is the number of nonnegative integer solutions of the equation
$m_1 + \cdots + m_\mu = m$.
In most examples of interest, only a few
types of compressed indecomposable dimension vectors
can contribute to a decomposition of $\dd$, and we can
obtain compact expressions for the number of orbits.
\\[1em]
{\bf Example.}  The dimension vector
$((n,n),(n,n),(1^{2n}))$, spherical of type ($\sp D_{2n+2}$).
The multiple flag variety consists of triples containing
two Lagrangian subspaces and one complete isotropic flag.
Let $c_n$ be the number of $\sp_{2n}$-orbits on the flag variety
$\cfl_{((n,n),(n,n),(1^{2n}))}$.

For $n=2$, there are only two families $(m_{\ee})$
with $\sum_{\ee} m_{\ee}\, \ee = ((22)(22)(1^4))$:\\
(1)\ $m_{((22)(22)(1^4))}=1$, all other $m_{\ee}=0$; \\
(2)\ $m_{((11)(11)(1001))}=m_{((11)(11)(0110))}=1$, all other $m_{\ee}=0$. \\
(Recall that a summand $\ee$ must be a triple of {\it symmetric} compositions.)
Since $\mu_{((22)(22)(1^4))}=2$ and
$\mu_{((11)(11)(1001))}=\mu_{((11)(11)(0110))}=\mu_{((11)(11)(11))}=5$, we get
$$
c_2 = \choos{2\!+\!1\!-\!1}{ 2\!-\!1} + \choos{5\!+\!1\!-\!1}{5\!-\!1}^2 =
27\ .
$$

For general $n$, any symplectic-indecomposable summand of
$((nn)(nn)(1^{2n}))$ must have compressed form:
$((11)(11)(11))$ with $\mu = 5$; or $((22)(22)(1^4))$ with $\mu = 2$.
A family $(m_{\ee})$ with $\sum_{\ee} m_{\ee}\, \ee = ((nn)(nn)(1^{2n}))$
has $m_{\ee} \leq 1$ and is equivalent to a
partition of the set $\{1,2,\ldots,n\}$ into subsets of sizes 1 and 2.
An easy computation gives the exponential generating function
$$
\sum_{n = 0}^{\infty} \frac{c_n x^n}{ n!} = 1 + 5x + \frac{27 }{ 2!} x^2 +
\cdots
= e^{x^2+5x}\ .
$$

\noindent
{\bf Example.} The dimension vector $((n,n), (1,2n-2,1), (1^{2n}))$,
spherical of type ($\sp E^{(b)}_{2n+3}$).
The multiple flag variety consists of triples containing
a Lagrangian subspace, a line, and a complete isotropic flag.

There are 3 compressed vectors that can serve as
symplectic-indecomposable summands:
$((11)(2)(11))$ with $\mu=2$;\
$((11)(11)(11))$ with $\mu=5$;\
and $((22)(121)(1^4))$ with $\mu=5$.
Furthermore, all the summands but one must be of the first kind.
The families $m_{\ee}$ are equivalent to choosing
either one or two elements in $\{1,2,\ldots,n\}$.
Thus the number of orbits is equal to
$$
n\!\cdot \choos{5\!+\!1\!-\!1}{5\!-\!1} \choos{2\!+\!1\!-\!1}{2\!-\!1}^{\!n-1}
+ \choos{n}{2}\!\cdot \choos{5\!+\!1\!-\!1}{ 5\!-\!1}
\choos{2\!+\!1\!-\!1}{ 2\!-\!1}^{\!n-2}
= 5 \cdot 2^{n-3} n (n+3) \ .
$$

\section{Proof of Theorems \ref{th:k=3} -- \ref{th:standard form}}

\subsection{Proof of Theorem \ref{th:GL-Sp relationship}(i)}

Although we discuss only alternating bilinear forms and the symplectic group,
the arguments in the proof of Theorem \ref{th:GL-Sp relationship} hold almost
verbatim for symmetric forms and the orthogonal group.

For a $k$-tuple of positive integers $(p_1,\ldots,p_k)$,
we define the semigroup $\Lam_{p_1,\ldots,p_k}$
and the categories $\FF_{p_1,\ldots,p_k}$ and $\CFF_{p_1,\ldots,p_k}$
analogously to their counterparts for $k=3$.
When there is no risk of ambiguity, we drop the subscripts
$(p_1,\ldots,p_k)$ and write $\Lam$, $\FF$, $\CFF$, etc.

Our first task is to show that each isomorphism class in $\CFF$
corresponds to a unique $\sp$-orbit
in a multiple  symplectic flag variety.
This is a consequence of the following general fact,
which is analogous to (but sharper than) a lemma of Richardson
\cite{richardson}, and which generalizes a result of Derksen
and Weyman \cite{dw}.

Consider a group $G$ acting
on a set $X$, and suppose we have
involutions $g \mapsto g^{\sigma}$ on $G$
and $x \mapsto \sigma x$ on $X$ such that for all $g$ and $x$,
$$
\sigma( g (\sigma x)) = (g^{\sigma}) (x).
$$
Let $G^{\sigma}$ and $X^{\sigma}$ denote the fixed point sets of $\sigma$.
Suppose further that:
\\[.2em]
(1) The group $G$ is a subgroup in $M^\times$, the
invertible elements of some finite-dimensional associative
algebra $M$ over an algebraically closed field $\kk$.
\\[.2em]
(2) The anti-involution $g\mapsto g^* := (g^{\sigma})^{-1}$
on $G$ extends to a $\kk$-linear anti-involution
$m \mapsto m^*$ of the algebra $M$.
\\[.2em]
(3) For any $x \in X^{\sigma}$,
the stabilizer $H = {\rm Stab}_G(x)$ is the group of
invertible elements of its linear span ${\rm Span}_{\kk}(H)$ in $M$.

\begin{proposition}
\label{pr:abstract orbit correspondence}
Let $G$ and $X$ be as above.
If two points in $X^\sigma$ are $G$-conjugate then they are
$G^\sigma$-conjugate.
\end{proposition}
{\it Proof.} As a preliminary, let us show that
if $H = {\rm Span}_{\kk}(H)^\times$ then
any element $h \in H$ with $h^* = h$ can be written as
$h = k^2=k^* k$ for some $k=k^* \in H$.
Since $M$ is finite-dimensional,
the subalgebra $\kk[h]\subset M$
is isomorphic to a quotient of a polynomial ring,
$\kk[h] \cong \kk[t]/(p(t))$, where
$p(t)$ is the minimal polynomial of $h$.
Furthermore, since $h$ is invertible, $p(0)\neq 0$.
Now, by the usual theory of finitely-generated $\kk[t]$-modules,
$h$ has a square root in
$\kk[h]$:  that is, $h = k^2$ for some $k = q(h) \in \kk[h]$.
Clearly $k^*=k$, and also
$k \in {\rm Span}_{\kk}(H)^{\times} = H$, as desired.

Now, let $x \in X^\sigma$, and let
$H$ be the stabilizer of $x$ in $G$.
Consider any point $g x \in X^{\sigma}$ for $g \in G$.
We have
$$
gx = \sigma(gx) = \sigma g (\sigma x) = g^{\sigma} x,
$$
and hence $g^* g x = x$.
Thus the element $h := g^* g$ lies in $H$ and satisfies $h^* = h$.

By the above discussion, we can write $h = k^* k$ for some $k \in H$.
Therefore we have
$
(g^{\sigma})^{-1} g = (k^{\sigma})^{-1} k,
$
and $gk^{-1} = (gk^{-1})^{\sigma} \in G^{\sigma}$.
But $k \in H = {\rm Stab}_G(x)$, so
$$
g x = gk^{-1} x \in G^{\sigma} x,
$$
and we are done.
\endproof

\medskip

The proof of Theorem \ref{th:GL-Sp relationship}(i) 
is now completed as follows.
Let $G = \gl(V)$, and $X$ be
a multiple $\gl$-flag variety $\fl_{\dd}(V)$
with symmetric dimension vector $\dd$.
Choose a symplectic form on $V$,
and let $m \mapsto m^*$ be the corresponding adjoint map on
$M = {\rm End} (V)$.
Take $\sigma: G \to G$ to be $\sigma (g) = {(g^*)}^{-1}$.
For $x \in X$, let $\sigma (x)$ be the multiple flag formed by
all orthogonals of subspaces in $x$.
Then $G^\sigma = \sp(V)$,
and $X^\sigma = \cfl_{\dd}$.
Conditions (1)-(3) of Proposition~\ref{pr:abstract orbit correspondence}
are clearly satisfied.
(In fact, (3) is valid for any multiple flag $x \in \fl_{\dd}$,
not necessarily $\sigma$-fixed,
since linear maps that preserve $x$ form a subalgebra in
$M = {\rm End} (V)$.)
By Proposition~\ref{pr:abstract orbit correspondence},
two elements of $\cfl_{\dd}(V)$ lie in the same $\gl(V)$-orbit
if and only if they lie in the same $\sp(V)$-orbit.
This is precisely what was to be shown.
\endproof

\subsection{Proof of Theorem \ref{th:GL-Sp relationship}(ii)}


We begin by recalling some well-known facts about $\FF$.
(See for example \cite[\S3.1]{mwz}.)
Although $\FF$ is not an abelian category,
it has a full faithful embedding into
the abelian category of \emph{quiver representations}, and
an indecomposable of $\FF$ is also indecomposable in the quiver category.
Therefore we can apply the Krull-Schmidt theorem
for abelian categories \cite{atiyah} to conclude that
an object in $\FF$ has a unique expression as a direct
sum of indecomposables.
Furthermore, the endomorphism ring of an indecomposable object
in $\FF$ is local:  that is, every endomorphism is either invertible
or nilpotent.


In what follows, we fix $J = (V;F)$, 
a multiple symplectic flag in $\CFF$.
Here $V$ is a symplectic vector space and
$F$ is a tuple of symplectic flags.  
We shall examine how $(V;F)$ decomposes
in the larger category $\FF$.

Any subspace $U\subset V$ induces a subobject $I=(U; F\cap U)$ in $\FF$.
A splitting of the vector space $V = U \oplus W$
induces a splitting of $(V;F)$ in $\FF$ if and only if
$F = (F\cap U) \oplus (F\cap W)$;  that is, for any space $A\subset V$
which is a member of any flag of $F$, we have $A = (A\cap U) \oplus (A\cap W)$.

We shall repeatedly use an elementary linear algebra fact.
For any vector space $V$ and subspaces $X,Y,Z \subset V$,
the following conditions are equivalent:
$$
(*) \ \
\begin{array}{ccc}
X\cap(Y\!+\!Z)= (X\!\cap\! Y)+(X\!\cap\! Z)
&\Longleftrightarrow&
X+(Y\!\cap\! Z)= (X\!+\!Y)\cap(X\!+\!Z)\\
& \Longleftrightarrow &
Z \cap (X\!+\!Y)= (Z\!\cap\! X)+(Z\!\cap\! Y)\ .
\end{array}
$$
If any of these conditions holds, the subspaces $X$, $Y$, and $Z$
generate a distributive lattice.

\begin{lemma}
\label{lem:splitting1}
Let $(V;F)$ be a symplectic multiple flag and $U\subset V$
a subspace such that $U \cap U^\perp = 0$ and $(U; F \cap U)$ has
symmetric dimension vector.
Then both multiple flags $(U; F \cap U)$ and $(U^\perp; F \cap U^\perp)$
are symplectic, and we have the splitting in $\CFF$:
$(V;F) = (U; F \cap U) \oplus (U^\perp; F \cap U^\perp)$.
\end{lemma}
\smallskip

\noindent {\it Proof.} 
The condition $U \cap U^\perp = 0$ means that 
the symplectic form on $V$ is non-degenerate on $U$	and on $U^\perp$.
Let $A, A^\perp \subset V$ be two orthogonal members of a flag in $F$.
Since $F$ has symmetric dimension vector, we have
$\dim(A\cap U) = \dim U -\dim(A^\perp \cap U)$.
Now, the orthogonal of $A^\perp \cap U$ in $U$ is
$$
(A^\perp\cap U)^\perp \cap U = (A+U^\perp)\cap U \supset A \cap U\ .
$$
Since the form is non-degenerate on $U$, the left and right sides have
the same dimension, and the containment is an equality:
$(A+U^\perp)\cap U = (A\cap U) + (U^\perp\cap U)$.  Thus by $(*)$,
$$
A = A \cap (U \oplus U^\perp) = (A\cap U) \oplus (A\cap U^\perp).
$$
This implies the splitting of $(V;F)$ in $\FF$.
In particular, $(U^\perp;F\cap U^\perp)$ also has symmetric dimension vector.
Thus the summands are symplectic. \endproof \smallskip


\begin{lemma}
\label{lem:splitting2}
Let $(V;F)$ be a symplectic multiple flag which splits in $\FF$ as
$(V;F)=(U;F\cap U) \oplus (W;F\cap W)$.  Then: \\[.3em]
(i)  We have the adjoint splitting
$(V;F) = (U^\perp; F\cap U^\perp) \oplus (W^\perp; F\cap W^\perp)$.\\[.3em]
(ii) The multiple flag $(W^\perp; F \cap W^\perp)$
is isomorphic to the dual $(U;F\cap U)^*$.
In particular, $(W^\perp; F \cap W^\perp)$ has dimension vector
opposite to that of $(U;F\cap U)$.
\\[.3em]
(iii) The projections $\alpha:V \cong U\oplus W \to U$ and
$\beta:V\cong U^\perp \oplus W^\perp \to W^\perp$ induce morphisms
$$
\alpha:(W^\perp; F\cap W^\perp) \to (U;F\cap U).
\quad \mbox{\rm and} \quad
\beta:(U;F\cap U) \to (W^\perp; F\cap W^\perp)
$$
\end{lemma}
\smallskip

\noindent {\it Proof.} (i) The space $V$ splits as $V = U^\perp \oplus
W^\perp$
because the symplectic form is non-degenerate.  Now, for orthogonal members
$A,A^\perp$ of $F$, we have
$A^\perp = (A^\perp\cap U)+(A^\perp\cap W)$, which is equivalent to
$A = (A+U^\perp)\cap(A+W^\perp)$.  But this is equivalent to
$A = (A\cap U^\perp)+(A\cap W^\perp)$ by $(*)$.
\\
(ii)  The symplectic form on $V$ induces a pairing between $U$ and $W^\perp$,
which is non-degenerate since $W^\perp \cap U^\perp = 0$.
The elements of $W^\perp$ which pair to zero with $A \cap U$ are exactly
$$
\begin{array}{rcl}
(A\cap U)^\perp \cap W^\perp
&=& (A^\perp + U^\perp)\cap W^\perp\\
&=& (A^\perp \cap W^\perp)+(U^\perp \cap W^\perp) = A^\perp \cap W^\perp\ .
\end{array}
$$
(Here we used (i) and $(*)$ for the second equality.)
Thus $W^\perp$ is isomorphic to $U^*$; and $A^\perp \cap W^\perp$, a member of
$F\cap W^\perp$, is isomorphic to $(U/A)^*$.
\\
(iii) For $A$ a member of $F$, we have
$$
\alpha(A\cap W^\perp) = ((A\cap W^\perp) + W)\cap U \subset (A+W)\cap U =
A\cap U\ .
$$
(Here we used $(*)$ for the last equality.)
Thus $\alpha(F\cap W^\perp)\subset F\cap U$,
and similarly $\beta(F\cap U)\subset F\cap W^\perp$.
\endproof  \smallskip

Now let $J = (V;F)$ be a symplectic multiple flag which is indecomposable in
$\CFF$,
but which splits in $\FF$ as $(V;F) = (U;F\cap U) \oplus (W;F\cap W),$
where $I = (U;F\cap U)$ is indecomposable in $\FF$.
We know that the endomorphism $\alpha\beta:(U;F\cap U)\to (U;F\cap U)$ 
defined in Lemma \ref{lem:splitting2} must be either invertible or nilpotent.
That is, at least one of $\alpha\beta$ and $\id_U-\alpha\beta$ is an
isomorphism.
\\[.2em]
{\bf Case 1.} The map $\alpha\beta$ is an isomorphism.
Since $\dim (U) = \dim (W^\perp)$, the maps $\alpha$ and $\beta$ are
invertible
and give isomorphisms between $(U;F\cap U)$ and
$(W^\perp; F \cap W^\perp)\cong (U;F\cap U)^*$.
In particular, $(U;F \cap U)$ has symmetric dimension vector.
We also have $U \cap U^\perp = {\rm Ker}(\alpha) = 0$.
Hence by Lemma \ref{lem:splitting1},  $(V;F)$ splits in $\CFF$
as $(V;F)=(U;F\cap U) \oplus (U^\perp; F\cap U^\perp)$.
But since $(V;F)$ is indecomposable in $\CFF$, we must have
$(V;F) = (U;F \cap U)$: that is, $(V;F)$ is $\FF$-indecomposable.
\\[.2em]
{\bf Case 2.} The map $\id_U - \alpha\beta$ is an isomorphism.
It follows from the definitions that
$$
{\rm Ker} (\id_U - \alpha \beta) =
U \cap (W^\perp + (U^\perp \cap W)) \ .
$$
Hence in our case, we have
$$
U \cap W^\perp \subset U \cap (W^\perp + (U^\perp \cap W)) = 0 \ ,
$$
so that the sum $U\oplus W^\perp$ is direct.

Furthermore, consider the projections $\gamma: U\oplus W^\perp \to U$ and
$\delta: U^\perp \oplus W^\perp \to U^\perp$.
We may easily check that $\gamma$ restricts to the following isomorphism
(with inverse $\delta$):
$$
\gamma: (U+W^\perp)\cap(U+W^\perp)^\perp 
= (U+W^\perp)\cap U^\perp \cap W
\stackrel{\sim}{\to} 
U \cap (W^\perp + (U^\perp \cap W)) = 0.
$$
Thus the symplectic form on $V$ is non-degenerate on $U+W^\perp$.

Finally, note that for any member $A$ of $F$, we have
$$
\begin{array}{rcl}
\dim( A\cap (U+W^\perp)) &\leq& \dim(U+W^\perp) - \dim(A^\perp\cap(U+W^\perp)) \\
&\leq& \dim U + \dim W^\perp -\dim (A^\perp\cap U) - \dim(A^\perp\cap W^\perp) \\
&=& \dim(A\cap U)+\dim(A\cap W^\perp)
\end{array}
$$
where the last equality is because 
$(W^\perp;F\cap W^\perp) \cong (U;F\cap U)^*$ by Lemma \ref{lem:splitting2}. 
Thus $A \cap(U+W^\perp) = (A\cap U) + (A\cap W^\perp)$, and we have
$$
(U+W^\perp; F\cap(U+W^\perp))\ \cong\ (U;F\cap U) \oplus (W^\perp; F\cap W^\perp),
$$
with symmetric dimension vector.

As before, we may now apply
Lemma \ref{lem:splitting1} to conclude that
$U \oplus W^\perp$ induces a symplectic summand of $(V;F)$, so that 
$V = U\oplus W^\perp$
and $(V;F) \cong (U;F\cap U)\oplus (U;F\cap U)^*$. That is, $J \cong I
\oplus I^*$.

This concludes the proof of Theorem \ref{th:GL-Sp relationship}. \endproof

\subsection{Proofs of Theorems~\ref{th:k=3} and \ref{th:finite-type-list}}
\label{sec:proofs-classification}

%

Given  dimension vectors $\dd, \dd' \in \Lam_{p_1,\ldots,p_k}$,
we say $\dd'$ is a \emph{summand} of $\dd$
if $\dd-\dd' \in \Lam_{p_1,\ldots,p_k}$.

\begin{proposition}
\label{cor:symplectic summands}
Let $\dd \in \Lam_{p_1,\ldots,p_k}$ be a symmetric dimension vector.\\
(i) If some symmetric summand of $\dd$ is not $\sp$-finite,
then $\dd$ is not $\sp$-finite. \\
(ii)  If there exist only finitely many $\CFF$-indecomposable classes
whose dimension is a summand of $\dd$, then $\dd$ is $\sp$-finite.
\end{proposition}


\noindent {\it Proof.} (i)
If a summand $\dd'$ is $\sp$-infinite, there
are infinitely many distinct $\CFF$-classes of dimension $\dd'$.
Taking a direct sum with any $\CFF$-class of dimension
$\dd-\dd'$ we obtain (by the unique decomposition in $\CFF$)
an infinite family of distinct $\CFF$-classes of dimension $\dd$,
and thus infinitely many $\sp$-orbits on $\cfl_{\dd}$.\\
(ii) The $\sp$-orbits on $\cfl_{\dd}$ are in bijection with
families of integers $(m_I)_I$ with $\sum_I m_I\dim(I) = \dd$,
where $I$ runs over all $\CFF$-indecomposables.  The hypothesis ensures
that there are only finitely many such families.
\endproof \smallskip

\noindent {\bf Proof of Theorem~\ref{th:k=3}.}
If $\dd=(\aa_1\ldots\aa_k)$ with $k\geq 4$ were $\sp$-finite,
then $(\aa_1,\aa_2,\aa_3,\aa_4)$ would be as well.
Thus by Proposition~\ref{cor:symplectic summands}(i), it is enough to show
that the summand $\ff_0=((1^2),(1^2),(1^2),(1^2))$
is \emph{not} $\sp$-finite.
But in this case we have $\cfl_{\ff_0} = \fl_{\ff_0}$
since every vector in a symplectic vector space is isotropic.
So our statement follows from Theorem~\ref{th:GL-Sp relationship}
and the well-known fact that $\ff_0$ is of infinite $\gl$-type.
\endproof
\medskip

\noindent
{\bf Proof of Theorem~\ref{th:finite-type-list}.}
First let us prove that all the cases in the theorem are indeed
$\sp$-finite.  By Theorem \ref{th:GL-Sp relationship},
we must show there are finitely many $\CFF$-isomorphism classes
in dimension $\dd$ of type ($\sp A$)---($\sp Y$).
For all the types except the last, this follows from
the $\gl$-classification in \cite[Theorem~2.2]{mwz}.
That is, the $\dd$ of types ($\sp A$)---($\sp E^{(b)}$)
are all $\gl$-finite, meaning
there are only finitely many $\FF$-isomorphism classes
in dimension $\dd$; but there must be even fewer
$\CFF$-isomorphism classes.

It remains to show that every symmetric dimension vector
of type ($\sp Y$) (that is, of the form $\dd = ((1,2n-2,1), \bb, \cc)$
with $\ell (\bb) = 3$) is $\sp$-finite.
We will use the criterion of Proposition \ref{cor:symplectic summands}(ii).
First, notice that any symmetric summand of $\dd$ with even 
total dimension is of the same
type ($\sp Y$) (or the type ($\sp A$) already dealt with).
Thus any indecomposable summand of the form
$J = I \oplus I^*$ must have $\dim(I)=\ee$ with $\ee+\ee^{\opp}$
of type $\sp Y$.  But then it is easily seen
from the classification of \cite{mwz}
that $\ee$ is $\gl$-finite thus producing only finitely many 
$\CFF$-indecomposables in dimension $\ee+\ee^{\opp}$. 

Thus it only remains to limit the $\CFF$-indecomposables which
are also indecomposable in $\FF$.
For this purpose, let $\Vert \aa \Vert^2 = a_1^2 + \cdots + a_p^2,$
and for $\dd = (\aa, \bb, \cc) \in \Lam$ with
$|\aa|= |\bb| = |\cc| = 2n$, define the \emph{Tits quadratic form} $Q$ by:
$$
\begin{array}{rcl}
Q(\dd) &=& \dim \gl_{2n} - \dim \fl_{\dd} \\
&=& \frac{1}{2}(\Vert \aa\Vert^2 + \Vert \bb\Vert^2 + \Vert \cc\Vert^2 
- (2n)^2).
\end{array}
$$
The fundamental result of V. Kac \cite{kac} states that
if $Q(\dd) > 1$, then there is no indecomposable isomorphism
class of $\FF$ in dimension $\dd$; and if $Q(\dd)=1$, then there
is at most one indecomposable.

Thus, it only remains to check the
symplectic finite type property for those $\dd$
of type ($\sp Y$) with $Q(\dd) \leq 0$.
The following lemma is an easy calculation with the Tits form.

\begin{lemma}
\label{lem:imaginary roots of finite type}
Consider a compressed symmetric dimension vector
of the form $\dd = ((1,2n-2,1), \bb, \cc)$ with $\ell (\bb) = 3$.\\
(i)  If $Q(\dd) \leq 0$, then $\dd$ must be one of the following:
$$
\dd_1 = ((1,2,1), (1,2,1), (1^4)), \quad
\dd_2 = ((1,4,1), (2^3), (1^6))\ ;
$$
in fact, $Q(\dd_1) = Q(\dd_2) = 0$.
\\[.2em]
(ii)  If $Q(\dd)=1$, then $\dd$ must be one of the following:
$$
\dd_1^+=((1,2,1),(1,2,1),(1,2,1)), \quad \dd_2^+=((1,4,1),(2^3),(1,1,2,1,1))
\ . $$
\end{lemma}


We finish off the type ($\sp Y$) with the following statement,
whose proof we postpone until \S\ref{sec:proof of standard form}.

\begin{lemma}
\label{lem:d1 and d2}
  For $\dd=\dd_1$ or $\dd_2$ as above,
there are exactly two $\FF$-indecomposable classes
of dimension $\dd$ which lie in $\CFF$.
\end{lemma}
\medskip

Next we show that any $\dd$ \emph{not} on the list of Theorem
\ref{th:finite-type-list} is $\sp$-infinite.
First we observe the following analogue of~\cite[Lemma~3.5]{mwz}.

\begin{lemma}
\label{lem:imaginary roots of infinite type}
Let $\dd$ be a compressed triple of symmetric compositions of the same
even number.
Then exactly one of the following holds:
\item {(i)} $\dd$ belongs to one of the types $(\sp A)$---$(\sp Y)$ in
Theorem~\ref{th:finite-type-list}.
\item {(ii)} $\dd$ has a summand whose compressed form is one
of the following five dimension vectors:
$$\ff_1= ((1^4)(1^4)(1^4)),\
\ff_2 = ((1,2,1), (1^4), (1^4)),\ \ff_3 = ((2,2), (1^4), (1^4))
\ ,$$
$$\ff_4 = ((3,3), (2^3), (1^6)), \quad \ff_5 = ((2^3), (2^3), (2^3))
\ .$$
\end{lemma}

\noindent
{\it Proof.}  Consider the tree of implications:
\\
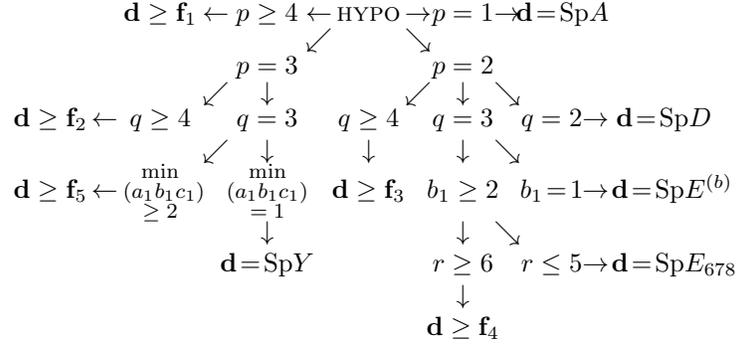
$$
\begin{array}{c@{\!}c@{\!}c@{\!}c@{\!}c@{\!}c@{\!}c@{\!}c@{\!}c@{\!}
c@{\!}c@{\!}l@{\!}l}
&& \dd\geq\ff_1 & \leftarrow & p \geq 4 & \leftarrow &\ \mbox{\sc hypo}\ &
\rightarrow & p=1 & \rightarrow & \dd\! =\! \sp A \!\!\!\!\!\\[-.2em]
&&&&& \swarrow && \searrow &&&&& \\[-.2em]
&&&& p=3 &&&& p=2 &&&& \\[-.2em]
&&& \swarrow & \downarrow &&& \swarrow & \downarrow & \searrow &&& \\[-.2em]
\dd \geq \ff_2\  & \leftarrow & q\geq 4 && q =3 &&
q\geq 4 && q=3 && q=2 &\rightarrow   \dd\! =\! \sp D\! \\[-.0em]
&&& \swarrow & \downarrow &&\downarrow && \downarrow & \searrow &&& \\[-.3em]
\dd \geq \ff_5\  & \leftarrow &
\mbox{\footnotesize $ \begin{array}{c}
{\rm min} \\[-.3em]\! (a_1b_1c_1)\!\!\! \\[-.3em] \geq 2
\end{array}$}
&&
\mbox{\footnotesize $ \begin{array}{c}
{\rm min} \\[-.3em]\!\!\! (a_1b_1c_1)\!\!\! \\[-.3em] =1
\end{array}$}
&&
 \dd \geq \ff_3 && b_1 \geq 2 && b_1\! =\! 1&\rightarrow \!\dd\! =\! \sp
E^{(b)}\!\\[-.3em]
&&&& \downarrow &&&& \downarrow & \searrow & && \\[-.0em]
&&&&\!\!\!\! \dd\! =\! \sp Y\!\!\!\! &&&& r\geq 6 && r\leq 5 &
 \rightarrow  \!\dd\!=\!\sp E_{678}\!\\
&&&&&&&&\downarrow&&&\\
&&&&&&&&\!\!\!\!\dd \geq \ff_4\!\!\!\!&&&
\end{array}
$$
\mbox{} \\[0em]
The root of the tree is our
\\[.2em]
{\sc hypothesis}: $\dd \in \Lam_{p,q,r}$
is a compressed triple of symmetric compositions of an even number, and
$1 \leq p \leq q \leq r$.
\\[.2em]
The arrows coming from a statement
point to all possible cases resulting from the statement.
We employ the abuse of notation $\dd = \sp A$, $\dd=\sp D$, etc to indicate
that $\dd$ belongs to the corresponding type
in Theorem~\ref{th:finite-type-list}.
We also write $\dd \geq \ff_i$ to indicate
that $\ff_i$ is a summand of $\dd$.
The lemma follows because every case ends in (i) or (ii),
and these conditions are clearly disjoint.
\endproof \smallskip

Now we can apply Proposition \ref{cor:symplectic summands}(i),
provided we show that each of the dimension vectors
$\ff_1, \dots, \ff_5$ is $\sp$-infinite.

For the first four of these cases, this is done by a simple dimension
count.  The symplectic group
and its isotropic Grassmannians have the following dimensions:
$$
\dim \sp_{2n} = n(2n+1), \quad \dim \cfl_{(k,2n-2k,k)} =
k(4n+1-3k)/2 \ .
$$
In particular, we have:
$$\dim \sp_4 = 10, \quad \dim \cfl_{(1^4)} = 4 $$
$$ \dim \cfl_{(1,2,1)} = \dim \cfl_{(2,2)} = 3,$$ \vspace{.3em}
$$\dim \sp_6 = 21, \quad \dim \cfl_{(1,4,1)} = 5, \quad
\dim \cfl_{(3,3)} = 6 \ ,$$
$$\dim \cfl_{(2^3)} = 7, \quad \dim \cfl_{(1^6)} = 9 \ .$$
We conclude that
$$\dim \cfl_{\ff_1} = 12 >
\dim \cfl_{\ff_2} = \dim \cfl_{\ff_3} = 11 >
\dim \sp_4 = 10\ ,$$
$$\dim \cfl_{\ff_4} = 22 > 21 = \dim \sp_6 \ .$$
Thus these four cases are of infinite symplectic type.

Now, $\dim \cfl_{\ff_5} = 21 = \dim \sp_6$, but we can easily see
that $\ff_5$ is also $\sp$-infinite by using
Theorem \ref{th:GL-Sp relationship}.
It is known (from counting dimensions)
that there are infinitely many non-isomorphic classes of
$\FF$ in dimension $\ff_5/2 = ((1^3),(1^3),(1^3))$.
Taking the direct sum of each such class with its dual,
we obtain infinitely many non-isomorphic classes of $\CFF$ in dimension
$\ff_5$,
meaning infinitely many $\sp_6$-orbits in $\cfl_{\ff_5}$.

This concludes the proof of Theorem~\ref{th:finite-type-list}.
\endproof

\subsection{Proof of Theorem \ref{th:standard form}}
\label{sec:proof of standard form}

To produce the list of finite-type symplectic indecomposables in Theorem
\ref{th:standard form}, we must consider all possible $\sp$-indecomposable
summands
for each type ($\sp A$)--($\sp Y$) in Theorem \ref{th:standard form}.
That is, we must refine our proof of
Theorem \ref{th:finite-type-list}, which
showed that there are only finitely many such indecomposables.
By Theorem \ref{th:GL-Sp relationship}, each $\sp$-indecomposable $J$ in
$\CFF_{pqr}$
is of the form $J = I$ or $J=I\oplus I^*$ for some $\gl$-indecomposable $I$
in $\FF_{pqr}$.

For all our types except ($\sp Y$), the corresponding dimension vectors are
not only
$\sp$-finite, but $\gl$-finite as well.
Thus we may use our list of $\gl$-indecomposables $I_{\dd}$ of finite type
from \cite[Thm 2.9]{mwz}.
Each $I_{\dd}$ is the unique $\gl$-indecomposable in its dimension.
(Notice that this implies $I_\dd^* \cong I_{\dd^\opp}$.)

Scanning the list of \cite{mwz}
we find that for every symmetric $\dd$ of type ($\sp A$)--($\sp E^{(b)}$),
the $\gl$-indecomposable $I_{\dd}$ lies in $\CFF$
(as demonstrated by the symplectic representatives given
in Theorem \ref{th:standard form}).
Therefore we can obtain all $\sp$-indecomposables for these types
either as:\\
(1)  $J=I_{\dd}$, where $\dd$ is symmetric
of type ($\sp A$)--($\sp E^{(b)}$) and $(I_{\dd})_\cpr$ is one of the
compressed $\gl$-indecomposables from \cite{mwz}; or \\
(2)  $J = I_{\ee}\oplus I_{\ee^\opp}$ where $\ee + \ee^\opp$ is of type
($\sp A$)--($\sp E^{(b)}$), $\ee$ is {\it not} symmetric, and
$(I_{\ee})_\cpr$, $(I_{\ee^\opp})_\cpr$ are among the compressed
$\gl$-indecomposables from \cite{mwz}.

Most of Theorem \ref{th:standard form} consists of
a systematic listing of these $J$ in compressed form.
(In compiling this list, one must remember that
even if $\dd=\ee+\ee^\opp$ is compressed, the summand $\ee$ might not be
compressed.)
\medskip

To complete our list, we must find the (compressed)
$\sp$-indecomposables of type ($\sp Y$).
Now, if $\dd$ is any dimension vector of type ($\sp Y$)
with summand $\ee+\ee^\opp$, then $\ee_\cpr$ belongs 
to one of the previous types
($\sp A$)--($\sp E^{(b)}$), and we can repeat the above procedure.

It remains to consider $\gl$-indecomposables $I$ of type ($\sp Y$)
which are symplectic.
By Kac's Theorem on indecomposables and
Lemma \ref{lem:imaginary roots of finite type},
such $I$ can occur only in the dimensions
$$
\dd_1 = ((1,2,1), (1,2,1), (1^4)), \quad
\dd_2 = ((1,4,1), (2^3), (1^6))\ ,
$$
\vspace{-1em}
$$
\dd_1^+=((1,2,1),(1,2,1),(1,2,1)), \quad
\dd_2^+=((1,4,1),(2^3),(1,1,2,1,1)).
$$
The bottom two dimensions have at most one indecomposable each.
We may easily check that the representatives $I_1^+ := I_{\dd_1^+}$ and
$I_2^+:=I_{\dd_2^+}$ given in Theorem \ref{th:standard form}
are $\gl$-indecomposable.
Indeed, the automorphisms of $I_1^+$ in $\gl_4$
are exactly the matrices of the form
$M(a,b)=$
{\tiny
$\left(\!\!\! \begin{array}{c@{\!}c@{\!}c@{\!}c}
a\ \, &0\ \, &0\ \, &0\\
0\ \, &a\ \, &b\ \, &0\\
0\ \, &0\ \, &a\ \, &0\\
0\ \, &0\ \, &0\ \, &a
\end{array}\!\!\! \right)
$ }\!
with $a\neq 0$.
Since this automorphism group contains no semi-simple elements other than
scalars,
$I_1^+$ is indecomposable (cf.~\cite[Lemma 2.3]{kac}).  
Similarly for $I_2^+$.
It is clear that $I_1^+, I_2^+$ lie in $\CFF$.

Finally we must find all the symplectic objects among the
(infinitely many) $\gl$-indecomposables in dimensions $\dd_1$ and
$\dd_2$.  This is the content of the lemma left unproved
in \S\ref{sec:proofs-classification}:
\\[.5em]
{\bf Proof of Lemma \ref{lem:d1 and d2}}.
Our technique is to fiber $\fl_{\dd_i}$ over $\fl_{\dd_i^+}$ by
dropping one subspace from the complete flag.
In other words, consider the {\it contraction functor}
$\pi_i:\FF_{p,q,r}\to \FF_{p,q,r-1}$ which sends an object $(V;A,B,C)$
to $(V;A,B,C')$, where $C'$ is obtained from $C$ by dropping the $i$th
subspace $C_i$.

\begin{lemma}
\label{lem:projection of indecomposable}
Suppose $(V;A,B,C)$ is an object of $\FF_{p,q,r}$ whose contraction
splits in $\FF_{p,q,r-1}$ as
$$
\pi_i(V;A,B,C) = \bigoplus_{j=1}^{c+1} (V_j;A^j,B^j,C^j),
$$
where $c = c_i+c_{i+1} = \dim(C_{i+1}/C_{i-1})$ and $V_1,\ldots, V_{c+1}
\subset V$.

Then $(V;A,B,C)$ itself splits in $\FF_{p,q,r}$ as
$$
(V;A,B,C)\ =\ (V_k; A^k, B^k, C^k)\, \oplus\,
(\tilde V_k; \tilde A^k, \tilde B^k, \tilde C^k)
$$
for some $k$, where $\tilde V_k = \bigoplus_{j\neq k}V_j$.
\end{lemma}

\noindent
{\it Proof.} It suffices to show
that for any subspaces $C' \subset C \subset C'' \subset V$
with $\dim(C''/C')=c$, if $C'$ and $C''$ split as
$$
C' = \bigoplus_{j=1}^{c+1} (C' \cap V_j) \qquad
C'' = \bigoplus_{j=1}^{c+1} (C'' \cap V_j),
$$
then $C$ splits as $C = (C\cap V_k) \oplus (C\cap \tilde V_k)$ for some $k$.

Since
$$c =  \dim  (C'' / C')  = \sum_{j = 1}^{c+1}
\dim  ((C'' \cap V_j)\, /\, (C' \cap V_j)) \ ,$$
there exists an index $k$ with $C'' \cap V_k = C \cap V_k= C' \cap V_k$.
Now, any $v \in C\subset C''$ has
a decomposition $v =v_k + \tilde v$ for $v_k \in C'' \cap V_k$,\,
$\tilde v \in  C'' \cap \tilde V_k$;
but then $v_k \in C\cap V_k$ and $\tilde v = v - v_k \in C \cap \tilde V_k$.
This implies the desired decomposition.
\endproof
\medskip

By the above lemma, any $\gl$-indecomposable $I$ of dimension $\dd_i$ can
split into
at most two summands when contracted to $\dd_i^+$; and if $I$ is symplectic,
then the contraction $I^+$ must be as well.
This leaves only a few possibilities for $I^+$.
We choose a representative for each possible $I^+$,
and insert an extra middle-dimensional subspace to lift it to dimension
$\dd_1$.
(This middle-dimensional space is automatically Lagrangian.)
In geometric terms, we consider the fibration
$\pi:\fl_{\dd_1}\to \fl_{\dd_1^+}$.
The automorphism group of $I^+$ acts on the fiber,
which is a projective line ${\bf P}^1$, and its orbits are the
$\CFF$-isomorphism classes of objects $I$ lying over $I^+$.  A given
$I$ is indecomposable exactly if its automorphism group contains no
semi-simple elements except scalars.

For $I$ of dimension $\dd_1$, the possible $I^+$ are:
$$
I_1^+,  \quad
I_{((110)(110)(110))}^\sym,   \quad
I_{((110)(110)(011))}^\sym,   \quad
I_{((110)(011)(110))}^\sym,   \quad
I_{((110)(011)(011))}^\sym \ .
$$
For $\dd_2$, they are:
$$
I_2^+,  \quad
I_{((120)(111)(11100))}^\sym,   \quad
I_{((120)(111)(10110))}^\sym,   $$\vspace{-1em} $$
I_{((120)(111)(01101))}^\sym,   \quad
I_{((120)(111)(00111))}^\sym \ .
$$
The analysis of these ten cases completes
the proof of Lemma \ref{lem:d1 and d2} and Theorem \ref{th:standard form}.
It turns out that only for $I^+=I_1^+$ and $I^+=I_2^+$ do we obtain
any $\gl$-indecomposable classes in the lifting:
two classes in each case.

We work out two typical cases. \\[.2em]
(i) $I^+= I_1^+$. The lifted objects $I=(V;A,B,C)$ are of the form:
$$V=(e_1,e_2,e_3,e_4)$$
$$
A= ((e_1+e_2+e_4)\subset (e_2,e_1+e_4,e_1+e_3)) \qquad
B = ((e_1)\subset (e_1,e_2,e_3)) ,
$$
$$
C = C(s:t)= ((e_4)\subset (se_2+te_3,e_4) \subset (e_2,e_3,e_4))\ .
$$
The automorphism $M(a,b)$ of $I_1^+$ takes $C(s,t)$ to
$C(s+\frac{b}{a} t:t)$.  Thus there are two orbits in the fiber: 
$I_{\dd_1}^1$ represented by $(s:t)=(1:0)$, 
and $I_{\dd_1}^2$ represented by $(s:t)=(0:1)$.
Both are indecomposable, since ${\rm Aut}(I^1_{\dd_1}) = {\rm Aut}(I_1^+)$,
and ${\rm Aut}(I^2_{\dd_1})$ consists of scalars. 
\\[.2em]
(ii) $I^+ = I_{((110)(110)(011))}^\sym$.  The lifted objects are:
$V = (e_1,e_2,e_3,e_4)$,
$$
A = ((e_1+e_2)\subset(e_1,e_2,e_3-e_4)),\qquad
B = ((e_1)\subset(e_1,e_2,e_3)),
$$
$$
C = C(s:t)= ((e_4)\subset(se_2+te_3,e_4)\subset(e_2,e_3,e_4))\ .
$$
Then
$$
{\rm Aut}(I^+)=\left\{ \left.
\mbox{ \small $
\left(\!\! \begin{array}{c@{\!}c@{\!}@{\!}c@{\!}c}
a\ &0\ \, &0\ &0\\[-.3em]
0\ &a\ \, &c\ &0\\[-.3em]
0\ &0\ \, &b\ &0\\[-.3em]
0\ &0\ \, &0\ &b
\end{array} \!\!\right)
$}
\ \right|\ a,b\neq 0 \, \right\}\ ,
$$
which has two orbits on the set of $I$, both clearly decomposable.
Thus, there are no indecomposables $I$ lying over $I^+$.
\smallskip

The other eight cases are similar.  \endproof

\end{document}